\documentclass[12pt]{amsart}

\usepackage{amssymb}

\usepackage{enumerate}

\usepackage{graphicx}

\makeatletter
\@namedef{subjclassname@2010}{%
  \textup{2010} Mathematics Subject Classification}
\makeatother

\newtheorem{thm}{Theorem}[section]
\newtheorem{cor}[thm]{Corollary}

\newtheorem{prop}[thm]{Proposition}

\theoremstyle{definition}
\newtheorem{defin}[thm]{Definition}
\newtheorem{rem}[thm]{Remark}
\newtheorem{exa}[thm]{Example}

\numberwithin{equation}{section}

\begin{document}


\baselineskip=17pt



\title[]{S\lowercase{trongly topologically transitive, supermixing, and hypermixing maps on general topological spaces}}
\author[M\lowercase{ahin} A\lowercase{nsari and} M\lowercase{ohammad} A\lowercase{nsari}]{M\lowercase{ahin} A\lowercase{nsari and} M\lowercase{ohammad} A\lowercase{nsari}}

\date{}

\begin{abstract} We give some basic properties of strongly topologically transitive, supermixing, and hypermixing maps on general topological spaces. 
Then we present some other results for which our mappings need to be continuous.
\end{abstract}


\maketitle

\section{\bf Introduction}
Let $X$ be a topological space. A map $f: X\to X$ is said to be
{\it hypercyclic} (or orbit-transitive) if there is some $x\in X$ for which
$$\textnormal{orb}(x,f)=\{f^nx: n=0, 1, 2. \cdots\}$$ is dense in
$X$.  Here $f^0=id_X$, the identity map on $X$, and $f^n=f\circ
f^{n-1}$ ($n\ge 1$). Thus, it is clear that non-separable topological spaces cannot support hypercyclic maps.
 If $\overline{\textnormal{orb}(x,f)}=X$ then $x$ is called a {\it hypercyclic point} for
$f$. The set of all hypercyclic points for $f$ is
denoted by $HC(f)$ and it is easy to see that, if $X$ has no
isolated point and $f$ is hypercyclic, then $HC(f)$ is dense in
$X$. If $HC(f)=X$ then $f$ is called {\it hypertransitive} (or minimal).
\par For a map $f: X\to X$, a set $\emptyset\neq E\subsetneq X$ is said to be $f$-invariant (or invariant under $f$) whenever $f(E)\subseteq E$. 
It is easy to verify that if $f$ is hypertransitive then $f$ lacks closed invariant subsets. The converse is true if 
$f$ is continuous.
\par A map $f:X\to X$ is called {\it topologically transitive} if, for
each pair of nonempty open sets $U, V\subseteq X$, there is some
$n\ge 0$ such that $f^n(U)\cap V\neq \emptyset$. One can
easily verify that if $X$ has no isolated point then any
hypercyclic map is topologically transitive. If, other than being empty of
isolated points, $X$ is a second countable Baire space then every
continuous topologically transitive map is hypercyclic (see Theorem $1.2$ of
\cite{bm} and its following remark).
 \par In \cite{ahkh}, the concept of strong topological transitivity has been introduced for 
 continuous linear operators on topological vector spaces. We can consider it for general maps on
 general topological spaces as well. A map $f:X\to X$ is called {\it strongly topologically transitive} if
$\bigcup_{n=0}^{\infty} f^n(U)=X$ for all nonempty open subsets
$U$ of $X$. Since a map $f:X\to X$ is topologically transitive if and only if
$\overline{\bigcup_{n=0}^{\infty} f^n(U)}=X$ for every nonempty open subset $U$ of $X$, we see that strong
topological transitivity implies topological transitivity.
\begin{rem}  In \cite[Definition $1$]{ka}, the finite union is considered to define
the notion of {\it strong transitivity}, i.e., a map $f: X\to X$ is called strongly transitive if, for 
any nonempty open set $U\subseteq X$, there exists a positive integer $s$ such that $\bigcup_{n=0}^s f^n(U)=X$.
 In this paper we only discuss strong topological transitivity (infinite union).
\end{rem}
A map $f$ is said to be {\it mixing} if, for
each pair of nonempty open subsets $U, V$ of $X$, there exists some
$N\ge 0$ such that $V\cap f^n(U)\neq \emptyset$ for all integers
$n\ge N$. Obviously, a mixing map is topologically transitive.
It is clear that, in the definition of a mixing map, there is no guarantee that $V$
intersects $A_N(U)=\bigcap_{n=N}^{\infty}f^n(U)$. Clearly, if $V$
intersects $A_N(U)$, then it also intersects each $T^n(U)$ ($n\ge
N$).
\par So, the property that, for any pair $U, V$ of
nonempty open subsets of $X$, there exists some $N\ge 0$ such that
$V\cap A_N(U)\neq \emptyset$ would be stronger than the mixing
property. It is equivalent to say that
$\bigcup_{i=0}^{\infty}\bigcap_{n=i}^{\infty}f^n(U)$ is dense in
$X$ for each nonempty open subset $U$ of $X$. 
\par Recently, the authors in \cite{a} have introduced and investigated
the notions of supermixing and hypermixing for continuous linear operators on topological vector spaces. 
We decided to consider these concepts for general maps supported
by general topological spaces.
\begin{defin} Let $X$ be a topological space. A map $f: X\to X$ is called
{\it supermixing} if, for each nonempty open
set $U\subseteq X$,
$$\overline{\bigcup_{i=0}^{\infty}\bigcap_{n=i}^{\infty}f^n(U)}=X.$$
 We say that $f$ is {\it hypermixing} if, for every nonempty
open subset $U$ of $X$, we have
$$\bigcup_{i=0}^{\infty}\bigcap_{n=i}^{\infty}f^n(U)=X.$$
\end{defin}
It is clear that any hypermixing map is supermixing and all
supermixing maps are mixing. Meanwhile, every hypermixing map is
obviously strongly topologically transitive.
\par The following example may be helpful when we want to think about the possible implications between the mentioned dynamical properties.
\begin{exa} (A) Let $X=\{a, b, c\}$ and $\tau=\{\emptyset, X, \{a,
b\}\}$ be a topology on $X$. \\
 ($\text{A}_1$) If $f:X\to X$ is the constant map $f(X)=\{c\}$
 then it is easy to see that $f$ is strongly topologically transitive but not
 mixing.\vspace*{.2cm}\\
 ($\text{A}_2$) Let $f$ be the function on $X$ for which $f(a)=f(b)=b$ and
 $f(c)=a$. Then $f$ is mixing but not strongly topologically transitive.\vspace*{.2cm}\\
($\text{A}_3$) If we define $f: X\to X$ by $f(a)=b, f(b)=c,
f(c)=a$, then it is not difficult to show that $f$ is mixing but
not
supermixing.\vspace*{.2cm}\\
($\text{A}_4$) Define the map $f$ on $X$ by $f(a)=a, f(b)=c,
f(c)=b$. Then it
can be easily seen that $f$ is supermixing but not hypermixing.\\
\\
(B) Let $X=\Bbb Z$ be equipped with the topology
$\tau=\{\emptyset, X, X\backslash\{0\}\}$ and the map $f: X\to X$
be defined by $f(m)=m+1$ ($m\in \Bbb Z$). One can readily verify
that $f$ is hypermixing.
\end{exa}
\par In Section $2$, we give some basic
results concerning strong topological transitivity, supermixing, and hypermixing properties. In Section $3$, we present some results involving other well-known 
dynamical properties for continuous maps.
\section{\bf Supermixing and hypermixing maps} 
The reader may have thought about the existence of a hypermixing
map on the topological space given in Example $1.3$ $(A)$, but,
while trying to find some, we felt that it is impossible. Then we
proved the following statement.
\begin{thm} Let $X$ be a topological space. If there is a nontrivial finite open set $U$ in $X$, then there
is no hypermixing map on $X$. In particular, if $X$ is a finite
set, then, equipped with any nontrivial topology, $X$ cannot
support a hypermixing map.
\begin{proof} Assume that $U$ is a nontrivial finite open subset of $X$.
If there is some $p\ge 1$ for which $U\subseteq f^p(U)$ then we
must have $U=f^p(U)$ since $U$ is a finite set. Then
$$\bigcup_{k=0}^{\infty}\bigcap_{n=k}^{\infty}f^n(U)=\bigcap_{n=0}^{p-1}f^n(U)\subseteq
U\neq X.$$ If for all $p\ge 1$, there is some $x(p)\in U\backslash
f^p(U)$ then $x(p)\notin \bigcap_{n=p}^{\infty} f^n(U)$. Since $U$
is a finite set, there must be some $x\in U$ and a strictly
increasing sequence of positive integers $(p_s)_s$ such that
$x=x(p_s)$ for all $s=1, 2, 3, \cdots$. Hence, $x\notin
\bigcap_{n=p_s}^{\infty}f^n(U)$ for all $s\ge 1$. Now, if $x\in
\bigcup_{k=0}^{\infty}\bigcap_{n=k}^{\infty}f^n(U)$ then there is
some $q\ge 0$ such that $x\in \bigcap_{n=q}^{\infty}f^n(U)$. On
the other hand, there is some $s_0\ge 1$ for which $p_{s_0}>q$.
Hence, $\bigcap_{n=q}^{\infty}f^n(U)\subseteq
\bigcap_{n=p_{s_0}}^{\infty}f^n(U)$ and so $x\in
\bigcap_{n=p_{s_0}}^{\infty}f^n(U)$, a contradiction.
\end{proof}
\end{thm}
A natural question which may often be asked whenever a new
property concerning a dynamical behavior of a map $f$ is
introduced, is that whether or not the maps $f^p$ ($p\ge 2$) too
have that property. The proof of the next result is exactly the same as that of Proposition $2.11$ of \cite{a}, but, 
for the sake of the reader's convenience, we bring it.
 \begin{prop} If $f:X\to X$ is a supermixing (hypermixing) map then $f^p$ is also supermixing
(hypermixing) for all integers $p\ge 2$.
\begin{proof} 
 Fix an integer $p\ge 2$. Suppose that $U$ is any nonempty open subset of $X$.
Then, for any $i\ge 0$, we have $$\{i, i+1, i+2,\cdots\}\supseteq
\{ip, (i+1)p, (i+2)p,\cdots\}.$$ Therefore,
$$\bigcap_{n=i}^{\infty}f^n(U)\subseteq \bigcap_{n=i}^{\infty}(f^p)^n(U),$$
for each $i\ge 0$. Hence,
$$\bigcup_{i=0}^{\infty}\bigcap_{n=i}^{\infty}f^n(U)\subseteq
\bigcup_{i=0}^{\infty}\bigcap_{n=i}^{\infty}(f^p)^n(U).$$
\end{proof}
\end{prop}
  The following result shows
that strongly topologically transitive maps
are nearly always onto.
\begin{prop} Let $X$ be a topological space in which there is some nonempty open set $U$ such that $\overline U\neq X$.
Then every strongly topologically transitive map on $X$ is onto.
\begin{proof} Suppose $f$ is a strongly topologically transitive map on $X$. Let $x\in X$.
If $x\notin U$ then $x\in f(X)$ since we have
$\bigcup_{n=0}^{\infty} f^n(U)=X$. Now, suppose that $x\in U$.
Since $\overline U\neq X$ there must be some nonempty open set $V$
in $X$ such that $U\cap V=\emptyset$. Thus, $x\notin V$ and the
equality $\bigcup_{n=0}^{\infty} f^n(V)=X$ shows that $x\in f(X)$.
\end{proof}
\end{prop}
If $f$ is a strongly topologically transitive {\it open} map then we can omit the
condition of having a nondense nonempty open set in the above
proposition. Indeed, for a given $x\in X$ and a fixed nonempty
open set $U$, we replace $U$ by the open set $f(U)$ in the
equality $\bigcup_{n=0}^{\infty} f^n(U)=X$ to obtain
$\bigcup_{n=1}^{\infty} f^n(U)=X$. Then it is clear that $f(X)=X$.
Briefly, we can give the following result.
\begin{prop} Every strongly topologically transitive open map is onto.
\end{prop}
Regrading the definitions of hypermixing and supermixing maps and the fact that $\bigcup_{i=0}^{\infty}\bigcap_{n=i}^{\infty}f^n(X)\subseteq f(X)$, 
the following result looks obvious.
\begin{prop} Every hypermixing map is onto and all supermixing maps have dense range.
\end{prop}
The following result says that nearly all of the hypermixing maps fail to be one-to-one.
\begin{prop} Let $X$ be a topological space in which there is some nonempty open set $U$ such that $\overline U\neq X$.
 Then no hypermixing map on $X$ is one-to-one.
\begin{proof} Let $f$ be a hypermixing map on $X$ and $V=X\backslash \overline U$. Then we have 
$\bigcup_{i=0}^{\infty}\bigcap_{n=i}^{\infty} f^n(U)=X= \bigcup_{i=0}^{\infty}\bigcap_{n=i}^{\infty} f^n(V)$.
Pick a point $x\in X$. Then there is a positive integer $N$, some $u\in U$, and some $v\in V$ such that $f^N(u)=x=f^N(v)$. Thus, $f$ is not
one-to-one since $u\neq v$.
\end{proof}
\end{prop} 
 \section{\bf Continuous maps}
 By proposition $2.5$, every supermixing map has dense range. For continuous maps on compact Hausdorff spaces, the range is the whole space.
  \begin{prop} Every continuous supermixing map on a compact Hausdorff space is onto.
  \begin{proof} Let $f$ be a continuous supermixing map on a compact Hausdorff space $X$. Then, by Proposition $2.6$, $f$ has dense range. But $f(X)$ is compact and hence
  closed. Thus, $f(X)=X$.
  \end{proof}
  \end{prop}
We devote the rest of our note to give some connections between the dynamical properties of continuous maps.
\par  Fix a continuous map $f:X\to X$. For a point $x\in X$, the
set $J^{mix}_f(x)=J^{mix}(x)$ is defined by
$$J^{mix}(x)=\{y\in X: \exists (x_n)_n \; \text{in}\; X\;\text{such that}\; x_n\to
x \;\text{and}\; f^nx_n\to y\}.$$ It is known that $J^{mix}(x)$ is
a closed $f$-invariant set, and that $f$ is mixing if and only if
$J^{mix}(x)=X$ for all $x\in X$, or equivalently, for all $x$ in a
dense subset of $X$ \cite[Exercise 1.4.4]{gp}.
\par If, for a set $B\subseteq X$, we put $J^{mix}(B)=\bigcup_{x\in
B}J^{mix}(x)$, then we give the following necessary and sufficient
condition for a continuous hypercyclic map to be mixing.
\begin{prop} Let $X$ be a topological space without isolated points and $f:X\to X$
be a continuous hypercyclic map. Then the following are
equivalent.\\
\hspace*{.3cm} $(1)$ $f$ is mixing\\
\hspace*{.3cm} $(2)$ $J^{mix}(HC(f))\cap HC(f)\neq \emptyset$
\begin{proof} $(1)\Rightarrow (2)$.
Suppose $f$ is mixing. Then $J^{mix}(x)=X$ for all $x\in X$. Thus
$J^{mix}(HC(f))\cap
HC(f)=X\cap HC(f)=HC(f)\neq \emptyset$. \\
 \hspace*{.3cm} $(2)\Rightarrow(1)$. Let $y\in J^{mix}(HC(f))\cap HC(f)$. Then
there is a hypercyclic point $x$ such that $y\in J^{mix}(x)$.
Since $J^{mix}(x)$ is a closed $f$-invariant set we have
$X=\overline{\text{orb}(y,f)}\subseteq
\overline{J^{mix}(x)}=J^{mix}(x)$. If we show that
$J^{mix}(f^kx)=X$ for all $k\ge 1$, then for all points $z$ in the
dense set $\text{orb}(x,f)$ we have $J^{mix}(z)=X$ which shows
that $f$ is mixing. To this end, fix $k\ge 1$. Since $y\in
J^{mix}(x)$, it is easy to see that $f^ky\in J^{mix}(f^kx)$. Then,
since $J^{mix}(f^kx)$ is $f$-invariant, we have
$X=\overline{\text{orb}(f^ky, f)}\subseteq
\overline{J^{mix}(f^kx)}=J^{mix}(f^kx)$ and we are done.
\end{proof}
\end{prop}
\begin{cor} Let $X$ be a topological space without isolated points and $f:X\to X$
be a continuous hypertransitive map. Then the following are
equivalent.\\
\hspace*{.3cm} $(1)$ $f$ is mixing\\
\hspace*{.3cm} $(2)$ $J^{mix}(X)\neq \emptyset$\\
\hspace*{.3cm} $(3)$ $\overline{J^{mix}(X)}=X$.
\begin{proof} $(1)\Rightarrow (2)$. If $f$ is mixing then, by Proposition $3.1$,
$$\emptyset\neq J^{mix}(HC(f))\cap HC(f)=J^{mix}(X)\cap
X=J^{mix}(X).$$
 \hspace*{.3cm} $(2)\Rightarrow (3)$. Suppose
$J^{mix}(X)\neq \emptyset$. Since $\overline{J^{mix}(X)}$ is a
closed $f$-invariant set and, meanwhile, $f$ lacks 
closed invariant subsets (since $f$ is
hypertransitive), we must have $\overline{J^{mix}(X)}=X$.\\
\hspace*{.3cm} $(3)\Rightarrow (1)$. Assume that
$\overline{J^{mix}(X)}=X$. Then
$$J^{mix}(HC(f))\cap HC(f)=J^{mix}(X)\cap X\neq
\emptyset.$$ Hence, $f$ is mixing by Proposition $3.1$.
\end{proof}
\end{cor}
Recall from the introduction that a continuous map is hypertransitive if and only if it does not admit closed invariant subsets.
\begin{prop} Let $X$ be a topological space and $f:X\to X$ be a continuous hypertransitive map.
Then the following are equivalent. \\
\hspace*{.3cm} $(1)$ $f$ is supermixing \\
\hspace*{.3cm} $(2)$
$\bigcup_{k=0}^{\infty}\bigcap_{n=k}^{\infty}f^n(U)\neq
\emptyset$, for any nonempty open set $U\subseteq X$.
 \begin{proof} $(1)\Rightarrow (2)$. If $f$ is supermixing then
  $\overline{\bigcup_{k=0}^{\infty}\bigcap_{n=k}^{\infty}f^n(U)}=X$,
 for any nonempty open set
$U\subseteq X$. \\
\hspace*{.3cm} $(2)\Rightarrow (1)$. Let $U$ be a nonempty open
subset of $X$ and put $A_k=\bigcap\limits_{n=k}^{\infty}f^n(U)$.
 We need to show that $\overline{\bigcup_{k=0}^{\infty}A_k}=X$.
It is clear that $f(A_k)\subseteq A_{k+1}$ and meanwhile
$\bigcup_{k=0}^{\infty}A_k=\bigcup_{k=1}^{\infty}A_k$ (since
$A_0\subseteq A_1$). Hence,
$$f(\bigcup_{k=0}^{\infty}A_k)=\bigcup_{k=0}^{\infty}f(A_k)\subseteq\bigcup_{k=0}^{\infty}A_{k+1}=\bigcup_{k=0}^{\infty}A_k.$$
Since $f$ is continuous, we have
$f(\overline{\bigcup_{k=0}^{\infty}A_k})\subseteq
\overline{\bigcup_{k=0}^{\infty}A_k}$ which shows that
$\overline{\bigcup_{k=0}^{\infty}A_k}=X$ because $f$ is
hypertransitive (and $\overline{\bigcup_{k=0}^{\infty}A_k}\neq
\emptyset$ by our assumption).
 \end{proof}
 \end{prop}
 Before we finish this note, regarding Proposition $3.2$,
 Corollary $3.3$, and Proposition $3.4$, let us clarify that
 for a continuous map, hypercyclicity (resp. hypertransitivity) by itself does not
 imply the mixing (resp. supermixing) property. To this end, let $X=\{a, b, c\}$ be equipped with the topology $\tau=\{\emptyset,
 X, \{a\}, \{b, c\}\}$. Define $f:X\to X$ by $f(a)=b$, $f(b)=f(c)=a$.
 Then it is easy to verify that $f$ is a continuous
 hypertransitive (and hence hypercyclic) map which is not supermixing (since it is not mixing).
 \par If the curios reader is thinking about the possibility of the implication
 (supermixing)$\Rightarrow$ (hypertransitive) for continuous maps, this example would
 be helpful: let $X=\{a, b, c\}$ and $\tau_1=\{\emptyset,
 X, \{a\}, \{a,b\}\}$ be a topology on $X$. Define $f:X\to X$ by $f(a)=a$,
 $f(b)=f(c)=c$. Then $f$ is a continuous supermixing map which is
 not hypertransitive.
 
\vspace*{.5cm}
\hspace*{.35cm} Mahin Ansari\\
\hspace*{.35cm} Shiraz University (Graduated)\\
\hspace*{.35cm} Shiraz, Iran\\
\hspace*{.35cm} {\it E-mail address}: \email{\small
mahin.ansari@ymail.com}\\
\\
\hspace*{.35cm} Mohammad Ansari\\
\hspace*{.35cm} Department of Mathematics\\
\hspace*{.35cm} Azad University of Gachsaran\\
\hspace*{.35cm} Gachsaran, Iran\\
\hspace*{.35cm} {\it E-mail address}: \email{\small
ansari.moh@gmail.com}, \email{\small mo.ansari@iau.ac.ir}

 \end{document}